\documentclass[11pt]{article}
\usepackage{amssymb,latexsym,amsmath, amsthm}
\usepackage{amsfonts, graphics}
\def\R{\mathbb R}

\def\N{\mathbb N}

\newcommand{\prf}{{\bf Proof.\ }}
\newcommand{\prfe}{\hspace*{\fill} $\Box$

\smallskip \noindent}

\newtheorem*{theorem}{Theorem}
\newtheorem*{cor}{Corollary}
\theoremstyle{remark}
\newtheorem*{remark}{Concluding Remarks}

\begin{document}

\title{The asymptotic behavior of solutions to the repulsive
       $n$-body problem}
\author{Gerhard Rein\\
        Fakult\"at f\"ur Mathematik, Physik und Informatik\\
        Universit\"at Bayreuth\\
        D-95440 Bayreuth, Germany\\
        email: gerhard.rein@uni-bayreuth.de}

\maketitle

\begin{abstract} 
  The $n$-body problem with a purely repulsive Coulomb interaction
  is considered. It is shown that for large times $t$
  the distance between any two particles grows linearly in $t$.
  The trajectory of each
  particle is asymptotically a straight line with a fixed velocity which is
  different for different particles.
\end{abstract}

\maketitle

\section{Introduction}
Consider $n$ point-charges with spatial coordinates
$x_i(t)\in \R^3$, $i=1,\ldots,n$, depending on time $t\in \R$. The charges are
all of the same sign, and for simplicity all the particles
have the same charge and the same mass which is set equal to unity; the
results below remain true without this assumption,
as long as all the charges have the same sign. The particles interact
via a repulsive Coulomb force
so that the evolution of the system is governed by the usual
$n$-body equations, but with the repulsive sign, written here
in first order form:
\begin{equation} \label{nbodysystem}
\dot x_i = v_i,\qquad
\dot v_i =\sum_{j=1,\ldots,n,  j\neq i,} \frac{x_i - x_j}{|x_i - x_j|^3},
\qquad i=1,\ldots,n.
\end{equation}
The energy of the system, which is recalled in the next section,
is conserved and positive definite. Hence the $v_i$ remain bounded and
$x_i\neq x_j$ for $i\neq j$ as long as the solution exists,
and hence every solution
exists globally in time. It is natural to ask how the system behaves for
$t\to\infty$ (and for $t\to -\infty$), and one expects that
due to the repulsive interaction
the system spreads out in space.
Indeed, in \cite{lape} it is conjectured that the distance between any two
particles goes to infinity for $t\to\infty$.

We shall prove that there are two constants $c_1, c_2 >0$
which depend only on the
initial data such that for $t$ sufficiently large,
\begin{equation} \label{mainclaim1}
c_1 t \leq   |x_i(t) - x_j(t)| \leq c_2 t,\ i\neq j.
\end{equation}
This will follow from a relation between the potential energy of the system
and its kinetic energy relative to its asymptotic configuration,
cf.\ the theorem in the next section.
This relation further implies that there exist parameters
$x_i^\ast, v_i^\ast \in \R^3$, $i=1,\ldots,n$, such that 
\begin{equation} \label{mainclaim2}
x_i(t) = x_i^\ast + t\, v_i^\ast + \mathrm{O}(\ln t),\ t\to \infty, 
\end{equation}
and $v_i^\ast\neq v_j^\ast$ for $i\neq j$.

Our analysis is based on arguments which are completely analogous to
those used in \cite{illrein} for the %analysis %of the asymptotics of the
plasma-physics case of the Vlasov-Poisson system, which can be thought of
as the continuum, mean-field limit of the system \eqref{nbodysystem}
as $n\to\infty$, see also \cite{Perth,Rein07}.
We are not aware that these arguments, which are simple enough,
have previously been used in the present context. The repulsive
$n$-body problem, which can also be viewed as a special case of the
so-called charged $n$-body problem with the repulsive Coulomb interaction
dominating the Newtonian gravitational attraction,
and which is certainly less important and less intriguing
than the classical, gravitational $n$-body problem, appears in various
places in the literature. In addition to \cite{lape} we mention
\cite{alf1,alf2,ate,benvid,chav,chav2}.

\section{Results and proofs}

Throughout this section we consider a solution to the system
\eqref{nbodysystem} with initial data and $n\in \N$ arbitrary but fixed.
The kinetic and potential energies of the solution are defined as
\[
E_\mathrm{kin}(t) 
:=
\frac{1}{2} \sum_{1\leq i\leq n} |v_i(t)|^2, \qquad
E_\mathrm{pot}(t) 
:=
\frac{1}{2} \sum_{1\leq i\neq j \leq n}\frac{1}{|x_i(t)-x_j(t)|},
\]
and the total energy is conserved,
\[
E(t) := E_\mathrm{kin}(t) + E_\mathrm{pot}(t) = E(0).
\]
This fact is well known and easily checked.
It implies that as long as the solution exists, the velocities
of the particles remain bounded and their relative
distances remain bounded away from zero,
\begin{eqnarray} \label{velbound}
  |v_i(t)| \leq \sqrt{2\,E(0)},\ |x_i(t)-x_j(t)|\geq\frac{1}{\sqrt{2\,E(0)}},\
  1\leq i\neq j \leq n,
\end{eqnarray}
so in particular the solution exists globally in $t$, which is well
known as well.

In addition to the standard energy quantities we define
\begin{equation}\label{ekinrel}
    E_\mathrm{kin}^\mathrm{rel}(t)
    := \frac{1}{2}\sum_{1\leq i\leq n}\left|v_i(t) - \frac{x_i(t)}{t}\right|^2.
\end{equation}
We refer to this quantity as the
{\em relative kinetic energy}; notice that
according to \eqref{mainclaim2} the limit of $x_i(t)/t$ for $t\to\infty$
is the asymptotic velocity of the $i$th particle.
We obtain the following identities.
\begin{theorem}
  \begin{itemize}
  \item[(a)]
    For all $t\in \R$,
    \[
    \frac{d}{dt}
    \left[t^2 E_\mathrm{kin}^\mathrm{rel}(t) + t^2 E_\mathrm{pot}(t)\right]
    =t\,E_\mathrm{pot}(t).
    \]
  \item[(b)]
    For all $t\geq 1$,
    \[
    E_\mathrm{kin}^\mathrm{rel}(t) + E_\mathrm{pot}(t)
    = \frac{C}{t}
    - \frac{1}{t}\int_1^t E_\mathrm{kin}^\mathrm{rel}(s)\,ds
    \]
    with $C>0$ independent of $t$.
  \end{itemize}
\end{theorem}
\noindent\prf
For the proof we recall the quantity
\[
I(t):=\frac{1}{2}\sum_{1\leq i \leq n} |x_i(t)|^2,
\]
which up to the factor  $1/2$ is the moment of inertia and
has the well known property that
\begin{eqnarray*}
  \frac{d^2}{dt^2}I(t)
  &=&
  \frac{d}{dt}\sum_{1\leq i \leq n} x_i(t)\cdot v_i(t)\\
  &=&
  2 \, E_\mathrm{kin}(t) +
  \sum_{1\leq i\neq j \leq n} x_i(t)\cdot \frac{x_i(t)-x_j(t)}{|x_i(t)-x_j(t)|^3}\\
  &=&
  2 \, E_\mathrm{kin}(t) + E_\mathrm{pot}(t) =
  2 \, E(t) - E_\mathrm{pot}(t).
\end{eqnarray*}
Using this and conservation of energy the assertion in part (a) follows:
\begin{eqnarray*}
\frac{d}{dt}\left[t^2 E_\mathrm{kin}^\mathrm{rel}(t)
    + t^2 E_\mathrm{pot}(t)\right]
&=&
\frac{d}{dt}\left[ t^2 E(t) + I(t) - t \frac{d}{dt} I(t)\right]\\
&=&
2 t E(t) - t \,\left(2 E(t) - E_\mathrm{pot}(t)\right) = t\,E_\mathrm{pot}(t).
\end{eqnarray*}
If we abbreviate $g(t) = t^2 E_\mathrm{kin}^\mathrm{rel}(t)$
and $h(t) = g(t) + t^2 E_\mathrm{pot}(t)$, we can for $t\neq 0$
rewrite the identity in (a) in the form
\[
\frac{d}{dt} h(t) = \frac{1}{t} \,h(t) - \frac{1}{t} \,g(t).
\]
We read this as a linear, inhomogeneous ODE for the function $h$
which we solve by the variation of constants formula, taking
$t=1$ as initial time; any other initial time $t_0>0$ would do just as well:
\[
h(t) = t\, h(1) - t\, \int_1^t \frac{g(s)}{s^2} ds.
\]
If we recall what $g$ and $h$ stand for and divide by $t^2$ this
proves part (b). \prfe

We exploit this theorem as follows.
\begin{cor} \label{nbodyasymp}
  \begin{itemize}
  \item[(a)]
    For all $t\geq 1$, $E_\mathrm{pot}(t) \leq C/t$,
    and Eqn.~\eqref{mainclaim1} holds.
  \item[(b)]
    For $i=1,\ldots,n$, limiting velocities
    $v_i^\ast =\lim_{t\to \infty} v_i(t)$ exist, and for all $t\geq 1$,
    $|v_i(t) - v_i^\ast| \leq C/t$. %\frac{C}{t}.
    If $i\neq j$, then $v_i^\ast\neq v_j^\ast$. 
  \item[(c)]
    In addition, there exist points $x_i^\ast \in \R^3$, $i=1,\ldots,n$,
    such that Eqn.~\eqref{mainclaim2} holds.
  \end{itemize}
\end{cor}
\noindent\prf
The estimate for $E_\mathrm{pot}(t)$ follows directly from part (b) of
the theorem.
It implies the lower bound in
\eqref{mainclaim1}, while the upper bound is obvious from
the bound for the velocities in \eqref{velbound}.

The lower bound on $|x_i(t) - x_j(t)|$ for $i\neq j$ and
Eqn.~\eqref{nbodysystem} imply that
$|\dot v_i(t)| \leq C/t^2$. Hence for all $t' \geq t\geq 1$,
\begin{equation} \label{veldiff}
  |v_i(t')-v_i(t)| \leq \frac{C}{t}.
\end{equation}
The bound for the velocities in \eqref{velbound} implies that
for any sequence $t_k \to\infty$ the sequence $v_i(t_k)$
has a convergent subsequence, but due to \eqref{veldiff}
its limit does not depend on the subsequence or the chosen
sequence so that $v_i(t)$ has a limit. Taking $t' \to \infty$
in \eqref{veldiff} yields the error estimate in part (b).
Clearly,
\[
x_i(t) = x_i(1) + (t-1)\, v_i^\ast + \int_1^t (v_i(s) - v_i^\ast)\,ds,
\]
and since by the error estimate in (b) the integral is $\mathrm{O}(\ln t)$,
Eqn.~\eqref{mainclaim2} is established as well.
If we express $|x_i(t) - x_j(t)|$ using this formula and observe that
this difference grows linearly in $t$ if $i\neq j$, we can conclude that
$v_i^\ast\neq v_j^\ast$ for $i\neq j$, and the proof is complete.
\prfe

\noindent
\begin{remark}
    \begin{itemize}
      \item[(a)]
  The linear
  growth rate in \eqref{mainclaim1} is sharp, since it gives both
  a lower and an upper bound for $|x_i(t) - x_j(t)|$. It is conceivable
  that the error terms in part (b) of the corollary can be improved, but
  probably not by the rather innocent estimates used here.
  \item[(b)]
  Part (b) of the theorem implies that
  $E_\mathrm{kin}^\mathrm{rel}(t)$ decays at least like $1/t$,
  but this rate is not sharp, since by the same identity,
  \[
  \int_1^\infty E_\mathrm{kin}^\mathrm{rel}(s)\,ds < \infty;
  \]
  if not, the right hand side of the identity would become negative for
  large $t$, which is not possible. If we substitute the asymptotics
  from part (b) of the corollary we see that for large $t$,
  $E_\mathrm{kin}^\mathrm{rel}(t)\leq C \ln^2 t /t^2$,
  which is consistent with the convergence of the above integral.
\item[(c)]
  It is far from clear whether the asymptotic behavior which was obtained
  in \cite{illrein}
  by completely analogous methods for the Vlasov-Poisson system
  in the plasma physics case is optimal. The author conjectures
  that it is not, and the wish to understand this issue lead
  him to the $n$-body considerations which he reports here.
    \end{itemize}
    \end{remark}

\end{document}